\documentclass[11pt]{article}
\usepackage{amssymb}

\newcommand{\newsection}{ \setcounter{equation}{0} \section}
\def\nn{\nonumber }
\def\la{\lambda}
\def\Lxi{\Lambda_{\xi}}

\def\T{T_{orb}}

\def\pat{\partial}
\def\LA{{\cal Z}_{\xi}} 

\providecommand{\href}[2]{#2}

\begin{document}

\begin{titlepage}
\rightline{LMU-TPW 2000-10}
\rightline{MPI-PhT/2000-19}

\vspace{4em}
\begin{center}

{\Large{\bf Realization of the Three-dimensional Quantum Euclidean Space by Differential Operators}}

\vskip 3em

{{\bf
S.\ Schraml and J.\ Wess }}

\vskip 1em

Max-Planck-Institut f\"ur Physik\\
F\"ohringer Ring 6, D-80805 M\"unchen\\[1em]

Sektion Physik, Universit\"at M\"unchen\\
Theresienstr.\ 37, D-80333 M\"unchen

\end{center}

\vspace{2em}

\begin{abstract} 
The three-dimensional quantum Euclidean space is an example of a non-commutative space 
that is obtained from Euclidean space by $q$-deformation. Simultaneously, angular 
momentum is deformed to $so_q(3)$, it acts on the $q$-Euclidean space that becomes a 
$so_q(3)$-module algebra this way. In this paper it is shown, that this algebra can be 
realized by differential operators acting on $C^{\infty}$ functions on $\mathbb{R}^3$. 
On a factorspace of $C^{\infty}(\mathbb{R}^3)$ a scalar product can be defined that 
leads to a Hilbert space, such that the action of the 
differential operators is defined on a dense set in this Hilbert space and algebraically 
self-adjoint becomes self-adjoint for the linear operator in the Hilbert space. The 
self-adjoint coordinates have discrete eigenvalues, the spectrum can be considered as a 
$q$-lattice.
\end{abstract}

\vfill

\end{titlepage}\vskip.2cm

\newpage
\setcounter{page}{1}

\newsection{Introduction}

The algebra of the Euclidean quantum space $\mathbb{R}^3_q$ \cite{FRT,Wess} has been discussed in 
\cite{LWW,3DQES}. Its defining relations are:
\begin{eqnarray}
  \label{Xalg}
  &&X^3X^+-q^2X^+X^3=0\nn\\
  &&X^3X^--q^{-2}X^-X^3=0\\  
  &&X^-X^+-X^+X^-=\la X^3X^3\nn\\
  &&\overline{X^3}=X^3,\quad \overline{X^+}=-qX^-\nn\\
  &&\la=q-q^{-1},\quad q\in\mathbb{R},\quad q\ge 1.\nn
\end{eqnarray}
This space is a $so_q(3)$-module algebra. The whole set of relations can be found in \cite{LWW,3DQES}. 
The generators of the quantum Lie algebra 
$so_q(3)$ are interpretated as angular momentum operators.

In this paper we are going to show that the elements of $\mathbb{R}_q^3$ and $so_q(3)$
can be represented by differential operators in $\mathbb{R}^3$.
We use polar coordinates $(r,\theta,\varphi)$, $(\xi\equiv\cos\theta)$ and find 
\begin{eqnarray}
  \label{Xo}
  X^3&=&r\xi\nn\\
  X^+&=&-re^{i\varphi}\sqrt{\frac{1-q^{-2}\xi^2}{1+q^{-2}}}\;q^{-2\xi\frac{\pat}{\pat\xi}-1}\\
  X^-&=&+re^{-i\varphi}\sqrt{\frac{1-q^{2}\xi^2}{1+q^{2}}}\;q^{2\xi\frac{\pat}{\pat\xi}+1}\nn
\end{eqnarray}
for the coordinates. 
The generators of the $q$-deformed orbital angular momentum 
are represented as follows:
\begin{eqnarray}
  \label{To}
  \T^3&=&\frac{1}{\lambda}\left(1-q^{4i\frac{\pat}{\pat\varphi}}\right)\nn\\
  \T^+&=&\frac{e^{i\varphi}}{q\la\xi}\sqrt{1-q^{-2}\xi^2}\;q^{-2\xi\frac{\pat}{\pat\xi}}\nn\\
      && -\frac{1}{q\la\xi}\sqrt{1-q^2\xi^2e^{4i\frac{\pat}{\pat\varphi}}}\;e^{i\varphi}\\
  \T^-&=&\frac{qe^{-i\varphi}}{\la\xi}\sqrt{1-q^{2}\xi^2}\;q^{2\xi\frac{\pat}{\pat\xi}}\nn\\
      && -\frac{q}{\la\xi}\sqrt{1-q^{-2}\xi^2e^{4i\frac{\pat}{\pat\varphi}}}\;e^{-i\varphi}\nn
\end{eqnarray}
All these operators (\ref{Xo}), (\ref{To}) have the expected limit (\ref{Xlim}), (\ref{Tlim}) 
for $q\rightarrow 1$.

A similar result for $so_q(3)$ has been obtained in \cite{CZ}, where its  
generators have been constructed in terms of generators of $so(3)$.

In \cite{3DQES,Fiore} the representation theory of the algebra has been studied.
It was found that the representation is unique (apart from a 
scaling factor for the radius) if we demand that the conjugation
properties 
\begin{eqnarray}
\label{conjugation}
  \overline{X^3}=X^3,\quad \overline{X^+}=-q X^-\nn\\
  \overline{T^3}=T^3,\quad \overline{T^+}=q^{-2} T^-
\end{eqnarray}
are represented by the conjugation of linear operators in a 
Hilbert space and that the equal sign in (\ref{conjugation}) includes the
domain on which the linear operators are defined. This domain is
supposed to be dense in the Hilbert space.

The spectrum of the linear operator $X^3$ that was found in 
\cite{3DQES} does not agree with the spectrum of the differential 
operator in (\ref{Xo}) if we consider it as differential
operator in $L^2$. We cannot expect the differential operators of
(\ref{Xo}) and (\ref{To}) to have the desired conjugation properties
as linear operators on $L^2$.

To obtain the representation found in \cite{3DQES} we use the following
strategy: we consider the space of $C^{\infty}$ functions on which the 
differential operators act. This space is made to an algebra by a
convolutionary product. This algebra has an ideal that is left 
invariant under the action of the differential operators 
(\ref{Xo}) and (\ref{To}). We consider the factor  space of the 
$C^{\infty}$ algebra with respect to this ideal. On this factor space we can 
define a scalar pro\-duct making it a Hilbert space. This is the 
representation space where the operators (\ref{Xo}) and (\ref{To})
have the desired properties.

To achieve this we start from a basis in the $C^{\infty}$ space
where the elements are the product of functions of $r$, $\xi$ and $\varphi$.

The set of $C^{\infty}$ functions that vanish at $r=r_0q^{4M+2}$ 
for all $M\in\mathbb{Z}$ forms an ideal ${\cal I}_r^{r_0}$
under pointwise multiplication in the set ${\cal F}_r$ of all $C^{\infty}$ functions $f(r)$, $r\in\mathbb{R}_+$.  
$r_0$ is an arbitrary positive parameter, the
scaling factor mentioned above.
Since the differential operators (\ref{Xo}) and (\ref{To}) do not change 
the radius, it would be enough, to consider functions that vanish at one 
fixed $r$. But if one wants to introduce real momenta, one needs a scaling 
operator $\Lambda$, such that $\Lambda R=q^4 R\Lambda$ \cite{LWW}. Therefore 
we consider a whole $q$-lattice in radial direction. 

We introduce the factor space
\begin{equation}
  H_r^{r_0}=\frac{{\cal F}_r}{{\cal I}_r^{r_0}}.
\end{equation}
The scalar product that makes this space a Hilbert space, which we denote by 
${\cal H}_r^{r_0}$, is formulated with the Jackson integral:
\begin{equation}
  (g,f)=\sum_{M=-\infty}^{\infty}q^{4M}g^*(q^{4M+2}r_0)f(q^{4M+2}r_0).
\end{equation}

The eigenvectors of the multiplication operator $r$ with the eigenvalues $r_0q^{4M+2}$ form a basis in
this Hilbert space. We denote these vectors by $u_M$:
\begin{equation}
  ru_M=r_0q^{4M+2}u_M.
\end{equation}

For the set ${\cal F}_{\xi}$ of functions $f(\xi)$ we proceed similarly. The product
is again the pointwise product of the functions. 

The ideal is generated by the functions that vanish at $\xi=\pm q^{2m-1}$ 
for all $m\in\mathbb{Z},\; m\le 0$, 
we call it ${\cal I}_{\xi}$. The representation space is the factor space
\begin{equation}
  H_{\xi}=\frac{{\cal F}_{\xi}}{{\cal I}_{\xi}}.
\end{equation}
The scalar product that makes it a Hilbert space ${\cal H}_{\xi}$ is again defined
with the help of the Jackson integral:
\begin{equation}
  (\psi,\phi)=\sum_{\sigma=\pm 1}\sum_{m=-\infty}^{0}q^{2m}\psi^{\ast}(\sigma q^{2m-1})\phi(\sigma q^{2m-1}).
\end{equation}
The eigenfunctions of $\xi$ in this Hilbert space will be denoted by $\chi_{\pm m_t}$:
\begin{equation}
  \xi \chi_{\pm m_t}=\pm q^{2m_t-1}\chi_{\pm m_t}.
\end{equation}

For the set $\widetilde{\cal F}_{\varphi}$ of functions $f(\varphi)$ we define the product by the convolution:
\begin{equation}
  (\widetilde{fg})_m=\widetilde{f}_m\widetilde{g}_m,
\end{equation}
where $\sim$ denotes the Fourier transformation. This defines an algebra. The
functions for which $\widetilde{f}_m=0$ for $m<\underline{m}$ form an ideal. 
We construct the factor space with the scalar product
\begin{equation}
  (h,g)=\sum_{m=\underline{m}}^{\infty}\widetilde{h}_m^*\widetilde{g}_m.
\end{equation}
This will lead to the representation space when we allow $\underline{m}$ to depend on
$m_t$.

In this space we have the following basis:
\begin{eqnarray}
  \psi_{M,m_t,m}&=&u_M\chi_{m_t}e^{im\varphi}\\
  M&=&-\infty\ldots\infty\nn\\  m_t&=&-\infty\ldots 0\nn\\
  m&=&m_t\ldots\infty\nn
\end{eqnarray}
To define the scalar product for functions $\psi(r,\xi,\varphi)$ we fouriertransform with 
respect to $\varphi$ to obtain $\widetilde{\psi}_m(r,\xi)$.
\begin{eqnarray}
  (\phi,\psi)&=&\sum_{M=-\infty}^{\infty}\sum_{m_t=-\infty\atop \sigma=\pm 1}^0\sum_{m=m_t}^{\infty}q^{4M}q^{2m_t}\times\\
  &&\times\widetilde{\phi}^*_m(q^{4M+2}r_0,\sigma q^{2m_t-1})\widetilde{\psi}_m(q^{4M+2}r_0,\sigma q^{2m_t-1}).\nn
\end{eqnarray}

\newsection{The ${\mathbf X}$-Algebra}

Our aim is to represent the algebra (\ref{Xalg}) in terms of differential
operators acting on $\mathbb{R}^3$. We use polar coordinates $(r,\theta,\varphi)$ and 
$\xi=\cos\theta$. An operator that will play a major role in this attempt is:
\begin{equation}
  \LA=\frac{1}{2}\left(\xi\frac{\pat}{\pat\xi}+\frac{\pat}{\pat\xi}\xi\right)=\xi\frac{\pat}{\pat \xi}+\frac{1}{2}.
\end{equation}
It is defined in such a way that when acting on $L^2$-functions in the 
common domain of $\LA^{}$ and $\LA^*$,
\begin{equation}
  \LA^*=-\LA^{}
\end{equation}
holds. The property of $\LA$ that we will use frequently is:
\begin{eqnarray}
  &&[\LA,\xi]=\xi,\nn\\
  &&e^{\alpha\LA}\xi e^{-\alpha\LA}=e^{\alpha}\xi.
\end{eqnarray}
We now make an ansatz:
\begin{eqnarray}
  X^3&=&r\xi\nn\\
  X^-&=&Arf(\xi)e^{-2\alpha\LA},\quad f(0)=1,\\
  X^+&=&Brg(\xi)e^{2\beta\LA},\quad g(0)=1.\nn
\end{eqnarray}
From the first two equations of (\ref{Xalg}) follows:
\begin{equation}
  e^{-2\alpha}=e^{-2\beta}=q^2,
\end{equation}
and from the third equation:
\begin{equation}
  AB\left\{f(\xi)g(q^2\xi)-f(q^{-2}\xi)g(\xi)\right\}=\la\xi^2.
\end{equation}
With the definition
\begin{equation}
  \label{phidef}
  \phi(\xi)=ABf(\xi)g(q^2\xi)
\end{equation}
this equation becomes:
\begin{equation}
  \phi(\xi)-\phi(q^{-2}\xi)=\la\xi^2
\end{equation}
and has the solution:
\begin{equation}
  \label{phisol}
  \phi(\xi)=\phi(0)+\frac{q^3}{1+q^2}\xi^2.
\end{equation}

It is natural to identify the radius $r$ with the invariant length
in $\mathbb{R}^3_q$:
\begin{equation}
  r^2=R^2\equiv X^3X^3-qX^+X^--q^{-1}X^-X^+.
\end{equation}
This determines $\phi(0)$:
\begin{equation}
  \label{phinull}
  \phi(0)=-\frac{q}{1+q^2}.
\end{equation}

To obtain $f$ as well as $g$ from $\phi$ we have to use the conjugation
property
\begin{equation}
  \overline{X^+}=-qX^-,
\end{equation}
that leads to
\begin{equation}
  Bg(q^2\xi)=-q\overline{A}\overline{f}(\xi)
\end{equation}
and, as a consequence of the definition (\ref{phidef}) of $\phi$:
\begin{equation}
  \label{phi}
  \phi=-q|Af|^2.
\end{equation}
Introducing $\Lxi\equiv q^{2{\cal Z}_{\xi}}$, such that $\Lxi\xi=q^2\xi\Lxi$, and 
combining (\ref{phisol}), (\ref{phinull}) and (\ref{phi}) leads to 
the result:
\begin{eqnarray}
  \label{Xreal}
  X^3&=&r\xi,\nn\\
  X^+&=&-re^{i\varphi}\sqrt{\frac{1-q^{-2}\xi^2}{1+q^{-2}}}\;\Lxi^{-1},\nn\\
  X^-&=&re^{-i\varphi}\sqrt{\frac{1-q^2\xi^2}{1+q^2}}\;\Lxi.
\end{eqnarray}
We have found a representation of the $X$-algebra. In the limit
$q\rightarrow 1$ we obtain:
\begin{eqnarray}
\label{Xlim}
  X^3&\to & r\cos\theta\nn\\
  X^+&\to & -\frac{1}{\sqrt{2}} r\sin\theta \;e^{i\varphi}\\
  X^-&\to & +\frac{1}{\sqrt{2}} r\cos\theta \;e^{-i\varphi}.\nn
\end{eqnarray}

\newsection{The ${\mathbf t}$-algebra}

There is a homomorphism of the $T$-algebra into the
$X$-algebra \cite{Fiore,CFM,3DQES}:
\begin{eqnarray}
  t^+&=&-\frac{1}{\la q^3}\sqrt{1+q^2}X^+(X^3)^{-1}\nn\\
  t^-&=&\frac{q^2}{\la}\sqrt{1+q^2}X^-(X^3)^{-1}\\
  t^3&=&\frac{1}{\la}\left(1+R^2(X^3)^{-2}\right).\nn
\end{eqnarray}
With (\ref{Xreal}) this $t$-algebra can be represented
by differential operators:
\begin{eqnarray}
  \label{top}
  t^3&=&\frac{1}{\la}(1+\xi^{-2})\nn\\
  t^+&=&\frac{1}{\la}e^{i\varphi}\xi^{-1}\sqrt{1-q^{-2}\xi^2}\;\Lxi^{-1}\\
  t^-&=&\frac{1}{\la}e^{-i\varphi}\xi^{-1}\sqrt{1-q^{2}\xi^2}\;\Lxi.\nn
\end{eqnarray}
These are differential operators acting on $C^{\infty}$ functions. We cannot
expect that they are defined on a dense set in $L^2$ such that the conjugation properties
\begin{equation}
  \overline{t^+}=q^{-2}t^-,\quad \overline{t^3}=t^3
\end{equation}
hold for the differential operators when the conjugation is 
identified with the conjugation of linear operators in the Hilbert space 
$L^2$: 
\begin{equation}
  \label{topcon}
  (t^+)^*=q^{-2}t^-,\quad (t^3)^*=t^3.
\end{equation}
From \cite{3DQES} we actually know that this cannot
be the case because we found that such a representation of the $t$-algebra
is unique and leads to a spectrum of $t_3$ with the eigenvalues
\begin{equation}
  \frac{1}{\la}\left(1+q^2q^{-4m_t}\right),\quad m_t\le0.
\end{equation}
For $\xi$ this implies that only eigenvalues $q^{2m-1}, m\le 0$
are allowed. Clearly, the Hilbert space of square integrable 
functions is not the Hilbert space that would lead to 
such a spectrum.

The spectrum of $t^3$ suggests that we should consider a factor
space of the $C^{\infty}$ functions of the following type:

Consider the linear space of $C^{\infty}$ functions on the interval $0<\xi<1$
\begin{equation}
  {\cal F}_{\xi}=\{f(\xi)|f\in C^{\infty}((0,1))\}
\end{equation}
and the subspace generated by the functions
\begin{equation}
  {\cal I}_{\xi}=\{h\in C^{\infty}([0,1])|h(\xi_m)=0\mbox{ for } \xi_m=q^{2m-1}, m\le 0\}.
\end{equation}
Under pointwise multiplication these functions form an algebra,
which we also call ${\cal F}_{\xi}$ and ${\cal I}_{\xi}$, respectively. The algebra ${\cal I}_{\xi}$ is 
an ideal of ${\cal F}_{\xi}$ and we can define the factor space
\begin{equation}
  \label{xifactorspace}
  H_{\xi}\equiv \frac{{\cal F}_{\xi}}{{\cal I}_{\xi}}.
\end{equation}
On this factor space $t^3$ is defined and has the desired 
eigenvalues. Eigenvectors of $t^3$ with the
eigenvalue $\frac{1}{\la}\left(1+q^2q^{-4m_t}\right)$
we shall denote by $\chi_{m_t}$.

Next we show that the ideal ${\cal I}_{\xi}$ is left 
invariant by the action of $t^+$ and $t^-$
\begin{eqnarray}
  t^+f(\xi)&=&\frac{e^{i\varphi}}{\la}\frac{1}{\xi}\sqrt{1-q^{-2}\xi^2}\;\Lxi^{-1}f(\xi)\nn\\
  &=&\frac{e^{i\varphi}}{\la}\frac{1}{\xi}\sqrt{1-q^{-2}\xi^2}\;q^{-1} f(q^{-2}\xi)
\end{eqnarray}
For $f\in {\cal I}_{\xi}$ follows $t^+f\in{\cal I}_{\xi}$.

Analogous:
\begin{equation}
  \label{t-op}
  t^-f(\xi)=\frac{e^{-i\varphi}}{\la}\frac{1}{\xi}\sqrt{1-q^2\xi^2}\;qf(q^2\xi).
\end{equation}
$t^-$ shifts the points $\xi=q^{2m-1}$
to the points $\xi=q^{2(m+1)-1}$. In the
definition of ${\cal I}_{\xi}$ we have $m\le 0$, but $t^-f(\xi)|_{\xi=q^{-1}}=0$,
as can be seen from (\ref{t-op}). Thus ${\cal I}_{\xi}$ is invariant
under the action of $t^-$ as well. It follows that $t^3$, 
$t^+$ and $t^-$  are well defined on $H_{\xi}$.

We are now going to show that we can define
a scalar product on $H_{\xi}$ to get an Hilbert space ${\cal H}_{\xi}$, 
such that $(t^+)^*=q^{-2}t^-$:
\begin{equation}
  (\psi,\phi)=\sum_{m=-\infty}^0\psi^*(\xi_m)\phi(\xi_m)q^{2m},\quad \xi_m=q^{2m-1}.
\end{equation}
We compute:
\begin{eqnarray}
  (\psi,t^+\phi)&=&\sum_{m=-\infty}^0q^{2m-1}\psi^*(\xi_m)\frac{e^{i\varphi}}{\la}\frac{\sqrt{1-q^{-2}\xi^2_m}}{\xi_m}\;\phi(\xi_{m-1})\nn\\
  &=&\sum_{m=-\infty}^{-1}q^{2m+1}\psi^*(\xi_{m+1})\frac{e^{i\varphi}}{\la}\frac{\sqrt{1-q^{-2}\xi_{m+1}^2}}{\xi_{m+1}}\;\phi(\xi_m)\nn\\
  &=&\sum_{m=-\infty}^0\left(\frac{e^{-i\varphi}}{\la}\frac{\sqrt{1-q^2\xi^2}}{\xi}q^{2{\cal Z}_{\xi}}\psi\right)^*(\xi_m)\phi(\xi_m)q^{2m-2}\nn\\
  &=&\frac{1}{q^2}(t^-\psi,\phi).
\end{eqnarray}
We have changed the summation index and extended the sum to include
$m=0$ because the summand vanishes there.

There is a differential operator that commutes with the differential 
operators $\vec{t}$ of (\ref{top}):
\begin{eqnarray}
  \label{xihatdef}
  &&\hat{\xi}\equiv\xi q^{2i\frac{\pat}{\pat\varphi}}\\
  &&[\hat{\xi},\vec{t}\;]=0,\quad [\hat{\xi},R]=0,\quad [\hat{\xi},\vec{X}]=0.\nn
\end{eqnarray}

We shall use this operator to represent the $K$-algebra of
ref \cite{3DQES} in terms of differential operators.

\newsection{The ${\mathbf K}$-Algebra}

The elements of the $K$-algebra as they were defined in ref \cite{3DQES}
all commute with the $X$ and the $t$ algebra. The $K$-relations are
\begin{eqnarray}
  \label{Kalg}
  q^2 K^3 K^+-q^{-2} K^+ K^3&=&(q+q^{-1}) K^+ \nn\\
  -q^{-2}K^3 K^-+q^2 K^- K^3&=&(q+q^{-1}) K^- \\
  q^{-1} K^+ K^--q K^-K^+&=&K^3=\frac{1}{\la}(1-\tau_k)\nn
\end{eqnarray}
and
\begin{equation}
  \label{Kcon}
  \overline{K^3}=K^3,\quad\overline{K^+}=-q^{-2}K^-.
\end{equation}
The representation of the $K$-algebra that has to be used for 
orbital angular momentum has eigenvalues of $\tau_k$ of the 
form $-q^{-4m_k-2}$, $m_k\ge0$. This motivates the ansatz
\begin{equation}
  \tau_k=-\hat{\xi}^2.
\end{equation}
From 
\begin{equation}
  \hat{\xi}e^{i\varphi}=q^{-2}e^{i\varphi}\hat{\xi}
\end{equation}
follows that a promising ansatz for $K^+$ and $K^-$ is:
\begin{eqnarray}
  \label{Kopans}
  K^+&=&h(\hat{\xi})e^{i\varphi}\nn\\
  K^-&=&j(\hat{\xi})e^{-i\varphi}.
\end{eqnarray}
It satisfies the first two relations of (\ref{Kalg}), the third
one leads to a recursion formula for
\begin{equation}
  \label{Jdef}
  J(\hat{\xi})=h(\hat{\xi})j(q^2\hat{\xi}).
\end{equation}
We find:
\begin{equation}
  \label{Jrek}
  J(\hat{\xi})-q^2J(q^{-2}\hat{\xi})=\frac{q}{\lambda}(1+\hat{\xi}^2)
\end{equation}
with the solution
\begin{equation}
  J(\hat{\xi})=-\frac{1}{\la^2}\left\{1+\beta\hat{\xi}-q^2\hat{\xi}^2\right\}
\end{equation}
$\beta$ being a free parameter, not determined by (\ref{Jrek}).

From the conjugation property (\ref{Kcon}) follows
\begin{equation}
  j(\hat{\xi})=-q^2\overline{h}(q^{-2}\hat{\xi})
\end{equation}
if 
\begin{equation}
  \overline{\left(i\frac{\pat}{\pat\varphi}\right)}=i\frac{\pat}{\pat \varphi}.
\end{equation}
The parameter $\beta$ is determined to be zero by the
orbital angular momentum condition. This can be seen by a direct
calculation following all the steps outlined in \cite{3DQES}.

The result is:
\begin{eqnarray}
  \label{Kop}
  K^+&=&\frac{e^{i\vartheta}}{q^2-1}\sqrt{1-q^2\hat{\xi}^2}\;e^{i\varphi}\nn\\
  K^-&=&-\frac{q^2e^{-i\vartheta}}{q^2-1}\sqrt{1-q^{-2}\hat{\xi}^2}\;e^{-i\varphi}.
\end{eqnarray}

We now turn our attention to the representation space of the $K$-algebra.
The operator $\xi$ is represented on the factor space ${\cal H}_{\xi}$ defined in
(\ref{xifactorspace}). The eigenvectors of $t^3$ were denoted by $\chi_{m_t}$,
from 
(\ref{top}) we learn:
\begin{equation}
  \xi\chi_{m_t}=q^{2m_t-1}\chi_{m_t},\quad m_t\le 0.
\end{equation}
The eigenvectors of $\hat{\xi}$ (\ref{xihatdef}) will be of the form 
$\chi_{m_t}e^{im\varphi}$:
\begin{equation}
  \label{xihatac}
  \hat{\xi}\chi_{m_t}e^{im\varphi}=q^{2(m_t-m)-1}\chi_{m_t}e^{im\varphi}.
\end{equation}
These are the eigenfunctions of $\tau_k$:
\begin{equation}
  \tau_k \chi_{m_t}e^{im\varphi}=-q^{4(m_t-m)-2}\chi_{m_t}e^{im\varphi}.
\end{equation}

From \cite{3DQES} we know that the eigenvalues of $\tau_k$ are $-q^{-4m_k-2}, m_k\ge 0$.
It follows that
\begin{equation}
  m=m_t+m_k,\qquad m_k\ge 0,\quad m_t\le 0.
\end{equation}

For $m_t$ fixed we find the condition $m\ge m_t$. This is in agreement
with the expression of $K^-$ in (\ref{Kop}). The operator $K^-$ 
changes the eigenvalue of $K^3$ from $m_k$ to $m_k-1$. When applied
to the eigenvector $m_k=0$, it should give zero.
\begin{eqnarray}
  \label{KzuNull}
  K^-\chi_{m_t}e^{im_t\varphi}&=&-q^2\frac{e^{-i\vartheta}}{q^2-1}\sqrt{1-q^{-2}\hat{\xi}^2}\;\chi_{m_t}e^{i(m_t-1)\varphi}\nn\\
  &=&0.
\end{eqnarray}
To show this we use (\ref{Kop}) and the action of $\hat{\xi}$
(\ref{xihatac}). In this way we could have found the condition $m_k\ge0$. 

We see that for a given eigenvalue of $\xi$, the space of functions on which the 
$K$-algebra is represented is given by functions of $\varphi$
with a truncated Fourier transformation:
\begin{equation}
  \label{truncFou}
  g(\varphi)=\frac{1}{\sqrt{2\pi}}\sum_{m=m_t}^{\infty}c_me^{im\varphi}.
\end{equation}
This space of functions is invariant under the $K$-algebra, it was 
sufficient to show this for $K^-$ (Eqn (\ref{KzuNull})) because $K^+$
shifts the eigenvalue from $m$ to $m+1$, and $K^3$ does not
change the eigenvalue.

If we define the product of two functions as the convolution defined 
as product of the Fourier transformation we again have constructed 
an ideal by (\ref{truncFou}).

The factor space of the $C^{\infty}$ functions of $\varphi, (0\le\varphi\le 2\pi)$ with respect
to this ideal we call $\widetilde{H}^{m_t}_{\varphi}$. 
On this space a scalar product is defined:
\begin{equation}
  (h,g)=\sum_{m=m_t}^{\infty}\widetilde{h}_m^*\widetilde{g}^{}_m,
\end{equation}
where $\widetilde{h}$ and $\widetilde{g}$ stand for the Fourier transformation of $h$
and $g$. With this scalar product the conjugation property of the $K$-algebra 
(\ref{Kcon}) becomes
\begin{equation}
  (K^3)^*=K^3,\quad (K^+)^*=-q^{-2}K^-.
\end{equation}
This can easily be verified by a resummation and the use of (\ref{KzuNull}).

\newsection{Orbital angular momentum}

Orbital angular momentum has been defined in \cite{3DQES}:
\begin{eqnarray}
  \tau_{orb}&=&\tau_t\otimes\tau_k\nn\\
  \T^3&=&t^3\otimes 1+\tau_t\otimes K^3\nn\\
  \T^+&=&t^+\otimes 1+\sqrt{-\tau_t}\otimes K^+\\
  \T^-&=&t^-\otimes 1-\sqrt{-\tau_t}\otimes K^-\nn
\end{eqnarray}
For the differential operators this becomes:
\begin{eqnarray}
  \label{Torbop}
  \tau_{orb}&=&q^{4i\frac{\pat}{\pat \varphi}}\nn\\
  \T^3&=&\frac{1}{\la}\left(1-q^{4i\frac{\pat}{\pat\varphi}}\right)\nn\\
  \T^+&=&\frac{1}{\la}e^{i\varphi}\frac{1}{\xi}\sqrt{1-q^{-2}\xi^2}\;\Lxi^{-1}\nn\\
      &&+\frac{1}{\xi}\frac{e^{i\vartheta}}{q^2-1}\sqrt{1-q^2\hat{\xi}^2}\;e^{i\varphi}\\
  \T^-&=&\frac{1}{\la}e^{-i\varphi}\frac{1}{\xi}\sqrt{1-q^{2}\xi^2}\;\Lxi\nn\\
      &&+\frac{1}{\xi}\frac{e^{i\vartheta}}{q^2-1}q^2\sqrt{1-q^{-2}\hat{\xi}^2}\;e^{-i\varphi}\nn\\
  &&\hspace{-3em}{\cal Z}_{\xi}=\xi\frac{\pat}{\pat \xi},\quad\Lxi=q^{2{\cal Z}_{\xi}},\quad \hat{\xi}=\xi q^{2i\frac{\pat}{\pat \varphi}}\nn
\end{eqnarray}
The representation space with the proper conjugation properties
has been constructed above.

It remains to show that in the limit $q\rightarrow 1$ the operators
in (\ref{Torbop}) tend to the generators of angular momentum. We take 
$q=e^h$ and study the limit $h\rightarrow 0$. It is easy to see that:
\begin{equation}
  \T^3\to -2i\frac{\pat}{\pat \varphi}.
\end{equation}
For $\T^+$ the limit is more involved, as the two parts of $\T^+$ 
have no individual limit.
\begin{eqnarray}
  \T^+&=&\frac{1}{2h}\frac{e^{i\varphi}}{\xi}\sqrt{1-q^{-2}\xi^2}e^{-2h(\xi\frac{\pat}{\pat \xi}+\frac{1}{2})}\nn\\
      &&+e^{i\vartheta}\frac{1}{2h}\frac{e^{i\varphi}}{\xi}\sqrt{1-q^{-2}\xi^2e^{4ih\frac{\pat}{\pat \varphi}}}.
\end{eqnarray}
For $e^{i\vartheta}=-1$ the singular parts cancel and we obtain:
\begin{eqnarray}
\label{Tlim}
  \T^+&\to&e^{i\varphi}\sqrt{1-\xi^2}\left\{-\frac{\pat}{\pat\xi}+\frac{\xi}{1-q^2}i\frac{\pat}{\pat\varphi}\right\}\nn\\
  &\to&e^{i\varphi}\left\{\frac{\pat}{\pat\theta}+\cot\theta\; i\frac{\pat}{\pat\varphi}\right\}.
\end{eqnarray}

An analogous result is obtained for $\T^-$. It is interesting to note
that the phase $e^{i\vartheta}$ in the expression for the orbital angular momentum
has been determined to be $e^{i\vartheta}=-1$
by the requirement that the limit $q\to 1$ exists. 
The condition that the differential operators have the correct $q\to 1$ limit 
restricts the choise of the operators, without this condition more operators 
would satisfy the algebra. 
   


\end{document}